\documentclass{amsart}

\usepackage{amsmath,amssymb,amsthm}

\usepackage[bb=dsserif]{mathalpha}
\usepackage{bm}

\usepackage{graphicx,tikz}

\newtheorem*{thm}{Theorem}
\newtheorem{proposition}{Proposition}

\theoremstyle{definition}

\theoremstyle{remark}

\begin{document}

\title[]{The first eigenvector of a Distance Matrix\\ is nearly constant}
\thanks{S.S. is supported by the NSF (DMS-2123224) and the Alfred P. Sloan Foundation.}

\author[]{Stefan Steinerberger}
\address{Department of Mathematics, University of Washington, Seattle, WA 98195, USA}
\email{steinerb@uw.edu}

\begin{abstract} Let $x_1, \dots, x_n$ be points in a metric space and define the distance matrix
$D \in \mathbb{R}^{n \times n}$ by ${D}_{ij} = d(x_i, x_j)$. The Perron-Frobenius Theorem implies that 
there is an eigenvector $v \in \mathbb{R}^n_{}$ with non-negative entries associated to the largest eigenvalue. We
prove that this eigenvector is nearly constant in the sense that the inner product with the constant vector $\mathbb{1} \in \mathbb{R}^n$
is large
$$  \left\langle v, \mathbb{1} \right\rangle  \geq \frac{1}{\sqrt{2}} \cdot \| v\|_{\ell^2}  \cdot \|\mathbb{1} \|_{\ell^2}$$
and that each entry satisfies $v_i \geq \|v\|_{\ell^2}/\sqrt{4n}$. Both inequalities are sharp. \end{abstract}

\maketitle

\section{Introduction}

 This paper is the result of trying to understand a curious phenomenon concerning the graph distance matrix (Problem 1, introduced below). We give a relatively simple argument (Proposition 1) which explains, in a particular setting, why one might expect to observe such a phenomenon. Proposition 1 is a conditional result and it is not a priori clear why or how often (meaning for `how many' graphs) the condition should be satisfied. A surprisingly large number of times it is indeed satisfied for reasons that are not presently understood and we will refer to this as Problem 2. The main contribution of the paper, besides introducing these two problems, is to show that there is indeed a general result (see \S 1.2) in metric spaces $(X,d)$ which is of intrinsic interest and explains the second phenomenon `up to constants'.  For many graphs that constant factor seems to rather close to $\sqrt{2}$ (the largest it could be) and this is currently not understood: we refer to \S 1.3.

\subsection{The first problem.} We start by explaining the first problem which first arose implicitly in \cite{stein}: let $G=(V,E)$ be a finite, connected graph on $n$ vertices $V = \left\{v_1, \dots, v_n\right\}$ and let
$D \in \mathbb{R}^{n \times n}$ denote the graph distance matrix defined  by ${D}_{ij} = d(v_i, v_j)$ where $d(v_i, v_j)$ is the graph distance between 
the vertices $v_i, v_j$.
\begin{quote} \textbf{Problem 1.} The linear system of equations
$$ D x = \mathbb{1} = (1,1,\dots, 1) \in \mathbb{R}^n$$
tends to have a solution $x \in \mathbb{R}^n$ for `most' graphs. Why is that?
\end{quote}
 It is certainly the case that for `many' graphs, the matrix $D$ simply happens to be invertible and this case is not particularly interesting. However, when the graph $G$ exhibits additional structure and/or symmetries, then this is often reflected in the non-invertibility of $D$, however, even then the linear system seems to have a solution for a constant right-hand side. We can make this quantitative: the database of graphs in Mathematica 12 contains
7642 connected graphs with $3 \leq |V| \leq 50$ vertices. Among those, the matrix $D$ is not invertible for 3155 of these graphs (this large percentage is due to the fact that the database contains more `structured' than `generic' graphs). However, $ D x = \mathbb{1}$ has a solution in all but 5 cases.\\
Using the Monte Carlo
method, we estimate that for Erd\H{o}s-Renyi graphs $G(n,p)$ the likelihood of a linear system associated to $G(50,p)$ not having a solution is (uniformly in $p$) less than $1\%$ and
decays rapidly for larger $n$. The phenomenon arose in \cite{stein} where solutions of $Dx = \mathbb{1}$ are used to construct a notion of curvature on graphs. Naturally, the effectiveness of such a notion hinges on whether a solution of $Dx = \mathbb{1}$ typically exists.
Our starting point is a spectral existence result.

\begin{proposition} Suppose $D \in \mathbb{R}_{\geq 0}^{n \times n}$ has eigenvalues $\lambda_1 > 0 \geq \lambda_2 \geq \dots \geq \lambda_n$ and eigenvector $Dv = \lambda_1 v$. If 
$$  1 - \left\langle v, \frac{\mathbb{1}}{\sqrt{n}}\right\rangle^2 < \frac{|\lambda_2|}{\lambda_1 - \lambda_2},$$
then $Dx = \mathbb{1}$ has a solution.
\end{proposition}
 The right-hand side is positive but need not be large. Note that $D$ vanishes on the diagonal and thus the sum of the eigenvalues is 0. Given our assumption about the sign of the eigenvalues, this implies $\lambda_1 = |\lambda_2| + \dots + |\lambda_n|$. We would thus expect that, generically, $\lambda_1 \gg |\lambda_2|$ which makes the right-hand side rather small. Hence, the applicability of Proposition 1 depends on whether $\left\langle v, \mathbb{1}/\sqrt{n} \right\rangle$ is close to 1, whether the first eigenvector $v$ is `nearly constant'.

\subsection{The Result} If $D \in \mathbb{R}^{n \times n}_{}$ is a matrix with nonnegative entries, then the Perron-Frobenius theorem guarantees the existence of an eigenvector $v \in \mathbb{R}^n_{}$ with nonnegative entries associated to the largest eigenvector.
\begin{thm}
Let $x_1, \dots, x_n  \in X$ be points in a metric space (not all identical) and $D \in \mathbb{R}^{n \times n}$ denote
the matrix given by $D_{ij} = d(x_i, x_j)$. Its Perron Frobenius eigenvector $v \in \mathbb{R}^{n}_{\geq 0}$, normalized to unit length $\|v\|=1$, satisfies
$$ \min_{1 \leq i \leq n} v_i \geq \frac{1}{2\sqrt{n}} \qquad \emph{and} \qquad \left\langle v, \mathbb{1} \right\rangle  \geq \frac{\sqrt{n}}{\sqrt{2}}.$$
\end{thm}

\begin{center}
\begin{figure}[h!]
\begin{tikzpicture}
\node at (0,0) {\includegraphics[width=0.17\textwidth]{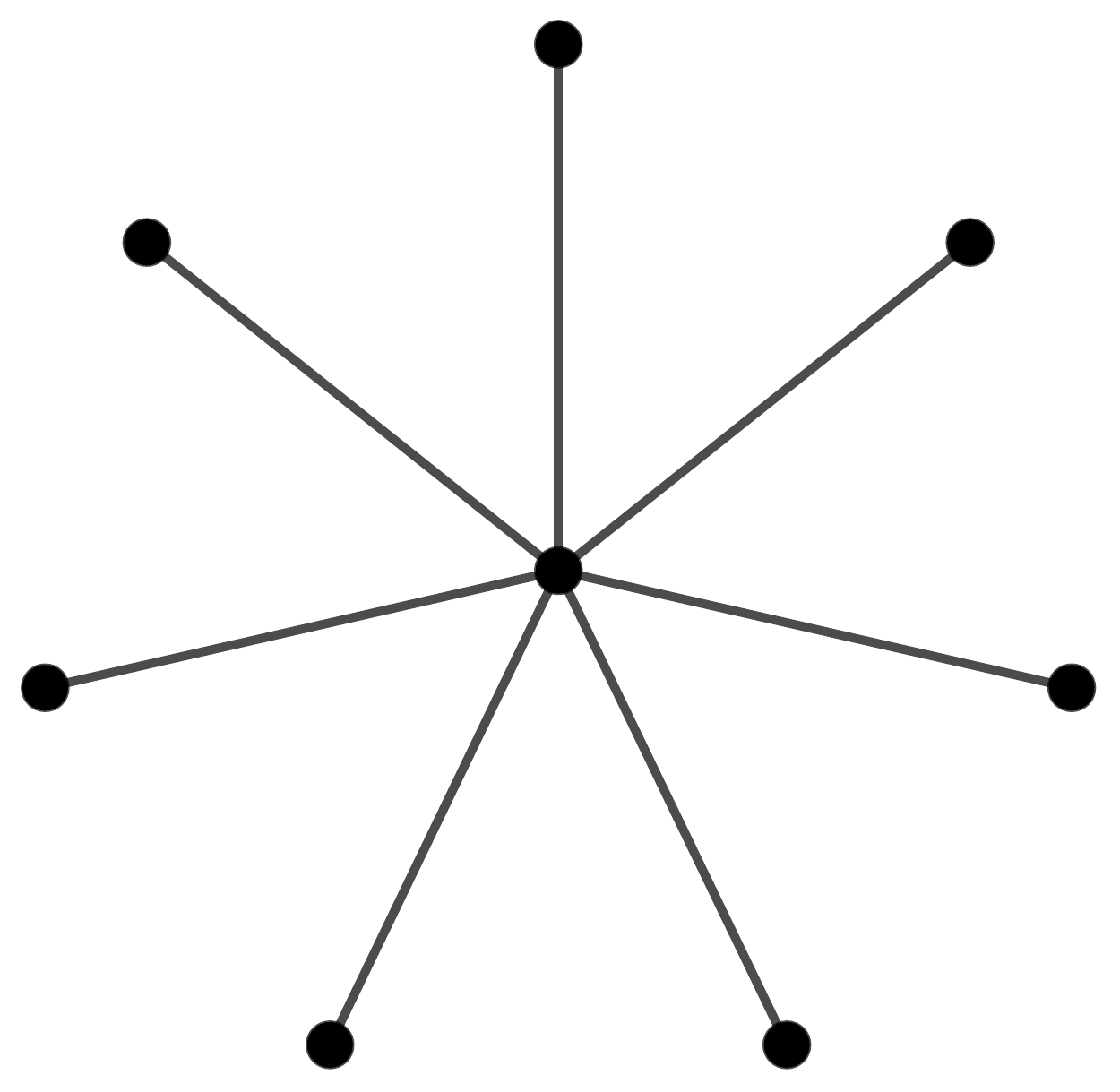}};
\filldraw (3,0) circle (0.05cm);
\filldraw (3.1,0.2) circle (0.05cm);
\filldraw (3.2,-0.1) circle (0.05cm);
\filldraw (3.3,0.2) circle (0.05cm);
\filldraw (3.2,0.25) circle (0.05cm);
\draw [<->] (2.9, -0.3) -- (3.4, -0.3);
\node at (3.15, -0.5) {$\ll \varepsilon$};
\filldraw (5,0) circle (0.05cm);
\draw [<->] (3.5, 0) -- (4.7, 0);
\node at (4.2, 0.2) {$1$};
\end{tikzpicture}
\caption{Sharpness: a star graph asymptotically for the first inequality, a metric space for the second inequality.}
\end{figure}
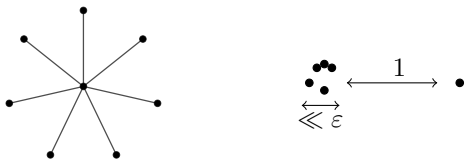
\end{center}
\vspace{-10pt}

We note that, by Cauchy-Schwarz, all $v \in \mathbb{R}^n_{\geq 0}$ of unit length $\|v\| = 1$ satisfy
$$ \max_{1 \leq i \leq n} v_i \geq \frac{1}{\sqrt{n}} \qquad \mbox{and} \qquad \left\langle v, \mathbb{1} \right\rangle  \leq \sqrt{n}$$
with equality if and only if all the entries of $v$ are $1/\sqrt{n}$. 

The Theorem can be interpreted as saying that the Cauchy-Schwarz inequality is close to being attained (which happens if the vectors $v$ and $\mathbb{1}$ are linearly dependent, meaning that $v$ is close to a constant vector). The two inequalities in the Theorem are sharp: if we consider the star graph, then the minimum of the largest eigenvector is attained
in the center and behaves asymptotically as described in the Theorem. 
For the second inequality we consider a metric space where $n-1$ points
are  distance $\leq \varepsilon$ (for suitable $\varepsilon$ depending on $n$) from each other with an additional point at distance $1$ from the first $n-1$ points. We refer to \S 2.3 for details.

\subsection{The second problem.} Returning to our original motivation, the Theorem by itself does not quite explain the phenomenon introduced above. The constant $1/\sqrt{2}$, while sharp, is not particularly close to 1. However, when returning to the graph setting, it does seem as if $\left\langle v, \mathbb{1}/\sqrt{n} \right\rangle$ is usually quite close to 1 and
often \textit{strikingly} so.  Among the 7642 connected graphs having $3 \leq n \leq 50$ vertices that are implemented in the Mathematica database, the \textit{average} value of $\left\langle v, \mathbb{1}/\sqrt{n} \right\rangle$ is 0.996.
Path graphs have a comparatively small constant that can be computed exactly (we refer to prior work of Ruzieh \& Powers \cite{ruz}). It is a nice fact (see \cite{ruz}) that the eigenvector behaves likes a hyperbolic cosine.
\begin{proposition}
Consider the path graph $P_n$ on $n$ vertices. Then
$$ \lim_{n \rightarrow \infty}  \left\langle v(P_n), \frac{\mathbb{1}}{\sqrt{n}} \right\rangle  = \frac{\sqrt{2} \sinh{c}}{\sqrt{c} \sqrt{c+ \cosh{c} \sinh{c}}} = 0.98261\dots ,$$
where $c$ is the positive number satisfying $c \tanh{c} = 1$.
\end{proposition}
A search in the Mathematica database identifies the `sun' graphs (also sometimes known as trampoline graph, see \cite{br}) as having the smallest constant within the library. A computation shows that
the limiting behavior of the constant for sun graphs is given by $(1/2 + 1/\sqrt{5})^{1/2} \sim 0.973$.
\vspace{-7pt}

\begin{center}
\begin{figure}[h!]
\begin{tikzpicture}
\node at (0,0) {\includegraphics[width=0.2\textwidth]{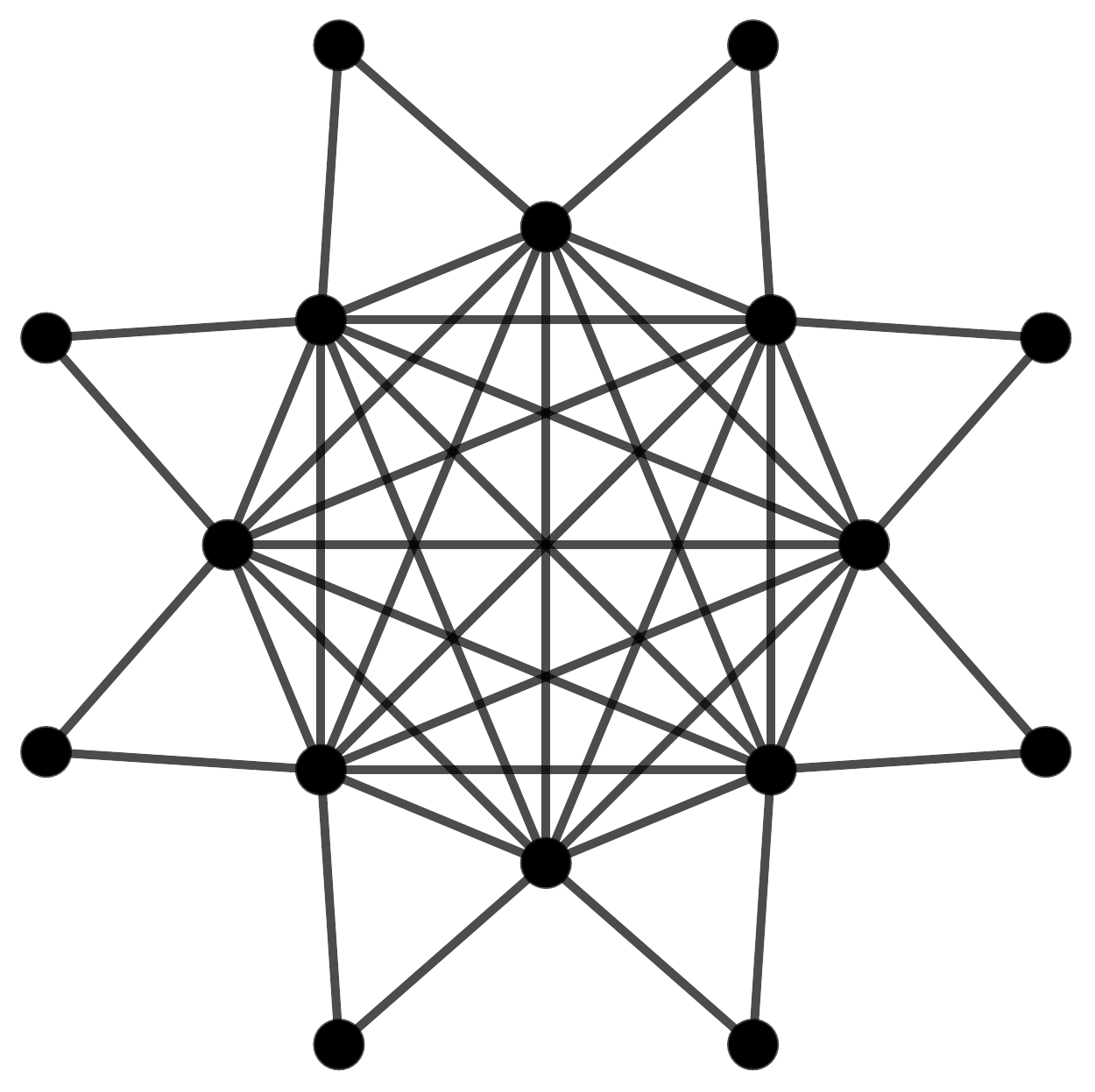}};
\node at (5.5,0) {\includegraphics[width=0.45\textwidth]{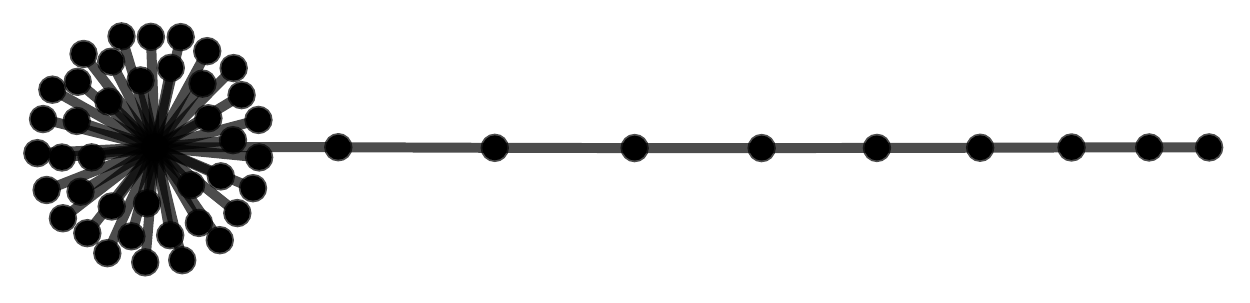}};
\end{tikzpicture}
\caption{The sun graph $S_8$ and a broom/comet.}
\end{figure}
\end{center}
\vspace{-13pt}
Motivated by the extremal metric space (see Fig. 1), it is natural to consider so-called \textit{broom graphs} or \textit{comets}  \cite{broom0, broom}. They are obtained by taking a star
and attaching a path to the central vertex. Numerics suggest that they might lead to a relatively small inner product $\left\langle v, \mathbb{1}/\sqrt{n} \right\rangle$ (certainly quite a bit smaller than the sun graphs and possibly even close to $1/\sqrt{2}$), however, their leading eigenvector does not seem to be known.  We conclude by formally stating the second problem.
\begin{quote}
\textbf{Problem 2.} If $G=(V,E)$ is an unweighted, connected graph on $n$ vertices and $v$ is the leading eigenvector of its distance matrix $D$, then it seem for most graphs $\left|\left\langle v, \mathbb{1} \right\rangle \right| \sim (1-\varepsilon) \sqrt{n}$ with a very small value of $\varepsilon$. Can this be made precise and/or quantified?
\end{quote}

\subsection{Related results}
 Let $(X,d)$ be a metric space and $x_1, \dots, x_n \in X$. The largest eigenvector of the distance matrix then describes the global maximum of the energy functional $E: \mathbb{R}^n \rightarrow \mathbb{R}$ given by
$$ E(a_1, \dots, a_n) = \frac{ \sum_{i,j = 1}^{n} d(x_i, x_j) a_i  a_j }{\sum_{i=1}^{n} a_i^2}.$$
This is naturally related to the study of energy integrals of the form
$$ E(\mu) = \int_{X} \int_{X} d(x,y) d\mu(x) d\mu(y),$$
where $\mu$ is a probability measure. The study of such measures is classical, we refer to the 1958 paper of Bj\"orck \cite{bj} as well as Alexander \cite{alex, alex2}, Alexander \& Stolarsky \cite{alex3}, Gross \cite{gross}, Wolf \cite{wolf2, wolf} and references therein. Both in the continuous \cite{bj} as well as in the discrete case \cite{steinembedding} the optimal measure $\mu$ tends to be localized in a rather small set contained in the boundary.
We observe that by switching from probability measures (essentially $L^1$) to $L^2$, the geometry of the problem changes substantially: having mass localized is actually rather expensive and it is beneficial to spread things out. Our main result makes this precise.
The distance matrix for graphs was introduced by Graham \& Pollak \cite{gra}.
Papendieck \& Recht \cite{pap} studied the largest possible entry of the principal eigenvector. 
Das \cite{po1, das2} considered the problem of bounding the extremal entries of the principal eigenvectors of the distance matrix.  The survey of Aouchiche \& Hansen \cite{survey} lists several more results about the principal eigenvalue and principal eigenvector, we also refer to the books of Stani\'c \cite{stanic}, Stevanovi\'c \cite{book}.
Distance matrices in the Euclidean setting were introduced by Schoenberg \cite{isaac} (see also Young \& Householder \cite{house}) and studied in \cite{ball, dis1, dis2, house}. We caution that the classical \textit{Euclidean Distance Matrix} (Menger \cite{menger}) is different and has $D_{ij} = \|x_i-x_j\|^2$. 

\section{Proofs}
\subsection{Proof of Proposition 1}
\begin{proof}
We consider the problem of maximizing the quadratic form $Q:\mathbb{R}^n \rightarrow \mathbb{R}$ 
$$ Q(v) =  \left\langle v, Dv \right\rangle$$
subject to the restriction
$$ \left\langle v, \mathbb{1} \right\rangle = 1.$$
Since $D$ has at least two positive entries, we conclude via $v= (1/n, \dots, 1/n)$ that 
$$ \sup_{\left\langle \mathbb{1}, v\right\rangle = 1} Q(v) > 0.$$
We will now argue that the maximum is attained. Using the symmetry of $D$ and the Spectral Theorem, we can decompose any vector into eigenvectors of $D$
$$ v =  \left\langle v, v_1 \right\rangle v_1 + \sum_{i=2}^{n} \left\langle v, v_i \right\rangle v_i.$$
Recalling that $\lambda_1 > 0 \geq \lambda_2 \geq \dots \geq \lambda_n$, we have
\begin{align*} Q(v) &=  \left\langle v, v_1 \right\rangle^2 \lambda_1 +  \sum_{i=2}^{n} \lambda_i \left\langle v, v_i \right\rangle^2 \leq   \left\langle v, v_1 \right\rangle^2 \lambda_1 + \lambda_2 \sum_{i=2}^{n} \left\langle v, v_i \right\rangle^2\\
&=   \left\langle v, v_1 \right\rangle^2 \lambda_1 + \lambda_2 \left( \|v\|^2 - \left\langle v, v_1 \right\rangle^2\right) = (\lambda_1 - \lambda_2) \left\langle v, v_1 \right\rangle^2  + \lambda_2 \|v\|^2.
\end{align*}
We can decompose the leading eigenvector $v_1$ (corresponding to $\lambda_1$) into its projection on $\mathbb{1}/\sqrt{n}$ and a remaining vector $ v_{1}^{*} \in \mathbb{R}^{n}$
$$ v_1 = \left\langle v_1, \frac{\mathbb{1}}{\sqrt{n}} \right\rangle \frac{\mathbb{1}}{\sqrt{n}} + v_{1}^{*}.$$
The Pythagorean theorem implies that $v_1^*$ has mean value 0 and thus
$$ 1 = \|v_1\|^2 = \left\langle v_1, \frac{\mathbb{1}}{\sqrt{n}} \right\rangle^2 + \|v_1^*\|^2.$$
Then
$$  \left\langle v, v_1 \right\rangle^2 = \left(   \left\langle v_1, \frac{\mathbb{1}}{\sqrt{n}} \right\rangle  \left\langle v, \frac{\mathbb{1}}{\sqrt{n}} \right\rangle +  \left\langle v,v_1^* \right\rangle \right)^2$$
and since the sum of entries of $v$ is $1$, this implies
$$   \left\langle v, \frac{\mathbb{1}}{\sqrt{n}} \right\rangle =   \left\langle v, \mathbb{1} \right\rangle \frac{1}{\sqrt{n}} = \frac{1}{\sqrt{n}}$$
and hence, for $v$ restricted to the hyperplane $x_1 + \dots + x_n=1$,
$$  \left\langle v, v_1 \right\rangle^2 =  \left(  \left\langle v_1, \frac{\mathbb{1}}{\sqrt{n}} \right\rangle \frac{1}{\sqrt{n}} + \left\langle v,v_1^* \right\rangle \right)^2.$$
Note that the first term is merely a constant, independent of $v$, and thus
$$  \left\langle v, v_1 \right\rangle^2 \leq  \left\langle v,v_1^* \right\rangle^2 + \mathcal{O}(\|v\| + 1).$$
 Using the bound from above and the Cauchy-Schwarz inequality,
 \begin{align*}
  Q(v) &\leq (\lambda_1 - \lambda_2) \left\langle v, v_1 \right\rangle^2  + \lambda_2 \|v\|^2 \\
  &\leq (\lambda_1 - \lambda_2) \left\langle v,v_1^* \right\rangle^2 + \lambda_2 \|v\|^2 + \mathcal{O}(\|v\| + 1) \\
  &\leq \left[(\lambda_1 - \lambda_2)\|v_1^* \|^2 + \lambda_2\right] \cdot \|v\|^2 + \mathcal{O}(\|v\| + 1).
  \end{align*}
 Therefore, if
 $$   \left( \lambda_1 - \lambda_2 \right) \|v_1^*\|^2 + \lambda_2 < 0,$$
 then there exists a global maximum that is attained because the quadratic form tends to $-\infty$ in all directions of the hyperplane $\left\langle v, \mathbb{1} \right\rangle =1$. Recalling that
 $$ \|v_1^*\|^2 = 1 -  \left\langle v_1, \frac{\mathbb{1}}{\sqrt{n}} \right\rangle^2$$
 we arrive at the desired condition. Let now $v_{\max}$ denote the location of such a maximum. It maximizes a function ($Q(v)$) subject to a constraint ($\left\langle v_{\max}, \mathbb{1} \right\rangle =1$)
 and therefore, for some Lagrange multiplier $\lambda \in \mathbb{R}$,
 $$  \nabla Q(v_{\max}) = 2Dv_{\max} = \lambda \cdot \mathbb{1}.$$
We note that since $\left\langle Dv_{\max}, v_{\max} \right\rangle = Q(v_{\max}) > 0$, we have $Dv_{\max} \neq 0$ and thus $\lambda \neq 0$. After rescaling, this implies
the existence of a solution of $Dx = \mathbb{1}$.
\end{proof}

\subsection{Proof of the Theorem}
\begin{proof} We start by noting that since not all the points $x_1, \dots, x_n \in X$ are identical, there exists at least a pair which has positive distance and thus the matrix $D$ has at least two positive entries (this is to ensure that we do not argue about the eigenvalues of the zero matrix). In particular, $\lambda_1 > 0$. Let $v=(v_1, \dots, v_n) \in \mathbb{R}^n$ denote the $\ell^2-$normalized Perron-Frobenius eigenvector all of whose entries are nonnegative. Then, for an arbitrary point $x_k$, we have 
\begin{align*}
 \lambda_1 = \left\langle Dv, v\right\rangle &= \sum_{i,j=1}^{n} d(x_i, x_j) v_i v_j \leq  \sum_{i,j=1}^{n} (d(x_i, x_k) + d(x_k, x_j)) v_i v_j \\
 &= \left( \sum_{j=1}^{n} v_j \right) \sum_{i=1}^{n} d(x_k, x_i) v_i +  \left( \sum_{i=1}^{n} v_i \right) \sum_{j=1}^{n} d(x_k, x_j) v_j \\
 &= 2  \left( \sum_{j=1}^{n} v_j \right) \sum_{i=1}^{n} d(x_k, x_i) v_i.
 \end{align*}
Using the $k-$th row of the equation $Dv = \lambda_1 v$, we note that
 $$  \sum_{i=1}^{n} d(x_k, x_i) v_i =  (Dv)_k = \lambda_1 v_k.$$
 Altogether, we conclude 
 $$  0 < \lambda_1 \leq 2  \|v\| _{\ell^1} \cdot \lambda_1 v_k = 2  \lambda_1  v_k \cdot  \left\langle v, \mathbb{1} \right\rangle$$
and thus, since $\lambda_1 > 0$,
  $$ v_k \geq \frac{1}{2 \left\langle v, \mathbb{1} \right\rangle}.$$
Then, using the Cauchy-Schwarz inequality,
 $$ v_k \geq \frac{1}{2 \left\langle v, \mathbb{1} \right\rangle} \geq \frac{1}{2 \|v\| \cdot \|\mathbb{1}\|} = \frac{1}{2\sqrt{n}}.$$
 This implies, in particular, $\left\langle v, \mathbb{1}/\sqrt{n} \right\rangle \geq 1/2$. However, one can do a little bit better: by noting that
$1 \leq k \leq n$ was arbitrary, we also obtain by summation over $k$ that
 $$ \left\langle v, \mathbb{1} \right\rangle \geq \frac{n}{2 \left\langle v, \mathbb{1} \right\rangle}$$
implying $ \left\langle v, \mathbb{1} \right\rangle^2 \geq n/2$ as desired.
\end{proof}

\textit{Remark.} The proof suggests that the only way for the argument to be close to sharp is if the triangle inequality is close to being an equation
for most cases. It seems interesting, albeit nontrivial, to see whether one can obtain quantitative improvements in settings where this
is not the case.

\subsection{Examples.} We will now discuss three relevant examples.\\
\textit{Star Graph.} The star graph $G_n$ is comprised of a central vertex $v$ which has an edge to $n$ other vertices. The symmetry of the graph is reflected in its leading eigenvector which we assume to have values $a$ in the center and $b$ in all other vertices.
 This leads to the problem of maximizing
$$ \frac{n(n-1) b^2 + n ab}{a^2 + nb^2} \rightarrow \max.$$
A computation shows that the optimal choice for $a$ is
$$ a = \frac{1}{2\sqrt{n}} +  \frac{5}{16} \frac{1}{n^{3/2}} + \dots$$
This is not too surprising: for the star graph, our argument is essentially optimal.\\
\textit{A metric space.}
As for the second inequality, let us suppose we have $n+1$ elements of which $n$ are all distance $\leq \varepsilon$ from each other with one element at distance $1$. The eigenvalue problem then consists of maximizing the expression
$$ \frac{\left\langle Dv, v \right\rangle}{\left\langle v, v \right\rangle } =  \frac{n a b + \mathcal{O}(n^2 a^2 \varepsilon)}{\sqrt{n a^2 + b^2}} \rightarrow \max,$$
where $a$ is the value in $n-1$ points that are distance $\leq \varepsilon$ from each other and $b$ is the value in the single point that is far away.
For $\varepsilon =0$, this is maximized for $a = 1/\sqrt{2n}$ and $b= 1/\sqrt{2}$. Basic stability theory implies that this persists for $\varepsilon$ sufficiently small. Thus
$$ \left\langle v, \frac{\mathbb{1}}{\sqrt{n+1}} \right\rangle = (1+o(1)) \cdot \frac{1}{\sqrt{2}}.$$
\begin{proof}[Path Graphs: Proof of Proposition 2] We use an exact expression for the first eigenvector of the distance matrix of $P_n$ due to 
Ruzieh \& Powers \cite{ruz}: they show that the entries $v_k$ for $1\leq k \leq n$ can be taken as
$$ v_k = \cosh{\left(\left(k - \frac{n+1}{2}\right) \theta \right)},$$
where $\theta$ is the positive solution of
$$ \tanh\left(\frac{\theta}{2} \right) \tanh\left(\frac{n \theta}{2} \right) = \frac{1}{n}.$$
A Taylor series expansion shows that for $x$ around 0
$$ \tanh{x} = x - \frac{x^3}{3} + \mathcal{O}(x^4)$$
and thus, as $n \rightarrow \infty$, we see that $\theta \sim n^{-1}$. Making the ansatz $\theta = (2c)/n$ for some unspecified constant $c$, we see that
\begin{align*}
 \frac{1}{n} &=  \tanh\left(\frac{\theta}{2} \right) \tanh\left(\frac{n \theta}{2} \right)  \\
 &=  \tanh\left(\frac{c}{n} \right) \tanh\left(c \right) = \frac{c}{n} \tanh\left(c \right) - \frac{c^3}{n^3} \tanh\left(c \right) + \dots
\end{align*}
Therefore, as $n \rightarrow \infty$, we have that the implicit constant $c_n$ in the ansatz $\theta_n = (2c_n)/n$ converges to the solution of the equation 
$c \cdot \tanh(c) = 1$. This solution is approximately $c \sim 1.2$. 
The limiting profile of the eigenvector is given by $\cosh{x}$ for $-c \leq x \leq c$ which, after some computation, implies the desired statement.
\end{proof}


\begin{thebibliography}{10}


\bibitem{alex} R. Alexander, Generalized sums of distances. Pacific J. Math. 56 (1975), no. 2, 297--304.

\bibitem{alex2} R. Alexander, Two notes on metric geometry.
Proc. Amer. Math. Soc. 64 (1977), p. 317--320.

\bibitem{alex3} R. Alexander and K. Stolarsky,
Extremal problems of distance geometry related to energy integrals.
Trans. Amer. Math. Soc. 193 (1974), 1--31.

\bibitem{survey} M. Aouchiche and P. Hansen, Distance spectra of graphs: a survey. Linear Algebra Appl. 458 (2014), 301--386. 

\bibitem{ball} K. Ball,
Eigenvalues of Euclidean distance matrices.
J. Approx. Theory 68 (1992), 74--82.


\bibitem{bj} G. Bj\"orck, Distributions of positive mass, which maximize a certain generalized energy integral, Ark. Mat., 3 (1958), 255--269.



\bibitem{dis1} E. Bogomolny, O.  Bohigas and C. Schmit,
Distance matrices and isometric embeddings. 
Zh. Mat. Fiz. Anal. Geom. 4 (2008), no. 1, 7--23, 202.

\bibitem{dis2} E. Bogomolny,  O. Bohigas and C. Schmit, 
Spectral properties of distance matrices. 
Random matrix theory.
J. Phys. A 36 (2003), no. 12, 3595--3616.

\bibitem{br} A. Brandst\"adt, A. V. B. Le and J. P.  Spinrad, Graph Classes: A Survey. Philadelphia, PA: SIAM, 1987.


\bibitem{po1} K. Das, 
Maximal and minimal entry in the principal eigenvector for the distance matrix of a graph. 
Discrete Math. 311 (2011), no. 22, 2593--2600.

\bibitem{das2} K. Das,
A sharp upper bound on the maximal entry in the principal eigenvector of symmetric nonnegative matrix
Linear Algebra Appl., 431 (2009), pp. 1340--1350

\bibitem{gra} R.L. Graham and H.O. Pollak. On the addressing problem for loop switching. Bell Syst. Tech. J.
50 (1971), 2495--2519.



\bibitem{gross} O. Gross, The rendezvous value of metric space. 1964 Advances in game theory pp. 49--53 Princeton Univ. Press, Princeton, N.J.

\bibitem{broom0} A. Ilic, On the extremal properties of the average eccentricity. Computers \& Mathematics with Applications, 64 (2012), 2877-2885.

\bibitem{menger} K. Menger, Untersuchungen uber allgemeine Metrik,
Math. Annalen 100, p. 75 -- 163 (1928)

\bibitem{pap} B. Papendieck and P. Recht, 
On maximal entries in the principal eigenvector of graphs. 
Linear Algebra Appl. 310 (2000), no. 1-3, 129--138.


\bibitem{ruz} S. Ruzieh and D.L. Powers, The distance spectrum of the path $P_n$ and the first distance eigenvector
of connected graphs, Linear Multilinear Algebra 28 (1990) 75--81.



\bibitem{isaac}  I. Schoenberg, Remarks to Maurice Frechet’s article Sur la definition axiomatique d’une classe
d’espace distancies vectoriellement applicable sur l’espace de Hilbert, Ann. of Math. 36 (1935), 724-732.

\bibitem{stanic} Z. Stani\'c,  Inequalities for graph eigenvalues. Vol. 423. Cambridge University Press, 2015.

\bibitem{stein} S. Steinerberger, Curvature on Graphs via Equilibrium Measures, arXiv:2202.01658

\bibitem{steinembedding} S. Steinerberger, Sums of Distances on Graphs and Embeddings into Euclidean Space, arXiv:2204.13278

\bibitem{book} D. Stevanovic, Spectral radius of graphs. Academic Press, 2014.

\bibitem{broom} D. Stevanovic and A. Ilic, Distance spectral radius of trees with fixed maximum degree. The Electronic Journal of Linear Algebra 20 (2010), p. 168-179.




\bibitem{wolf2} R. Wolf, On the average distance property and certain energy integrals.
Ark. Mat. 35 (1997), no. 2, 387--400.

\bibitem{wolf} R. Wolf, 
Averaging distances in real quasihypermetric Banach spaces of finite dimension. 
Israel J. Math. 110 (1999), 125--151.



\bibitem{house} G. Young and A. Householder, Discussion of a set of points in terms of their mutual distances, Psychometrika 3 (1938) 19--22.

\end{thebibliography}
\end{document}